\documentclass[12pt]{article}

\usepackage{amsmath}[1996/11/01]
\usepackage{amssymb,amsthm,amsfonts,latexsym,psfig,epsfig,here,graphics}
\usepackage{amscd, makeidx}

\def\biblitem#1{\bibitem{#1}}

\newcommand{\leftBra}{\{ \hspace*{-.10in} \{ }
\newcommand{\rightBra}{\} \hspace*{-.10in} \} }

\newcommand{\FF}{{\mathbb F}}
\newcommand{\beql}[1]{\begin{equation}\label{#1}}
\newcommand{\eeq}{\end{equation}}
\newcommand{\eqn}[1]{(\ref{#1})}
\numberwithin{equation}{section}
\newtheorem{thm}{Theorem}[section]

\theoremstyle{remark}
\newtheorem*{remark}{{\bf Remark}}

\newtheorem*{remarks}{{\bf Remarks}}

\theoremstyle{definition}
\newtheorem{eg}[thm]{Example}
\newtheorem{example}[thm]{Example}
\newtheorem{rems}[thm]{Remark}
\newtheorem{defn}[thm]{Definition}

\DeclareMathOperator{\Tr}{Tr}
\DeclareMathOperator{\tr}{tr}
\DeclareMathOperator{\ttt}{t}

\DeclareMathOperator{\Quad}{Quad}
\DeclareMathOperator{\End}{End}

\DeclareMathOperator{\Mat}{Mat}
\DeclareMathOperator{\rad}{rad}
\DeclareMathOperator{\Bil}{Bil}

\DeclareMathOperator{\Inv}{Inv}

\DeclareMathOperator{\Hom}{Hom}

\DeclareMathOperator{\GL}{GL}
\DeclareMathOperator{\Sp}{Sp}

\DeclareMathOperator{\cwe}{cwe}

\DeclareMathOperator{\diag}{diag}

\newcommand{\zentr}{\mbox{\large{\sf Y}\normalsize}}

\newcommand{\C}{\mathbb C}

\newcommand{\Q}{\mathbb Q}
\newcommand{\Z}{\mathbb Z}

\newcommand{\F}{\mathbb F}

\begin{document}

\LARGE
\begin{center}
{\bf Codes and Invariant Theory}
\end{center}
\normalsize

\begin{center}
G. Nebe, \\
Abteilung Reine Mathematik, \\
Universit\"{a}t Ulm, \\
89069 Ulm, Germany \\
Email: nebe@mathematik.uni-ulm.de \\
\vspace*{+.2in}

E. M. Rains, \\
Mathematics Department, \\
University of California Davis, \\
Davis, CA 95616, USA \\
Email: rains@math.ucdavis.edu \\
\vspace*{+.2in}

N. J. A. Sloane, \\
AT\&T Shannon Labs, \\
Florham Park, NJ 07932-0971, USA \\
Email: njas@research.att.com \\

\bigskip
January 1, 2003; revised October 6, 2003 \bigskip
\end{center}

\Large
\begin{center}
{\bf Abstract}
\end{center}
\normalsize
The main theorem in this paper is a
far-reaching generalization of Gleason's theorem
on the weight enumerators of codes which applies to 
arbitrary-genus weight enumerators of self-dual codes
defined over a large class of finite rings and modules.
The proof of the theorem uses a categorical approach,
and will be the subject of a forthcoming book.
However, the theorem can be stated and applied
without using category theory, and we illustrate it here 
by applying it to
generalized doubly-even codes over fields of characteristic 2,
doubly-even codes over $\Z/2^f\Z$,
and self-dual codes over
the noncommutative ring $\F_q + \F _q \, u$, where $u^2 = 0$..

\vspace{0.7\baselineskip}

\section{Introduction}

One of the most remarkable theorems in coding theory is
Gleason's 1971 theorem \cite{Gle70}
that the weight enumerator of a binary doubly-even
self-dual code is an element of the polynomial ring
generated by the weight enumerators of the Hamming code
of length $8$ and the Golay code of length $24$.
In the past thirty years many generalizations of this
theorem have been given that apply to other familes
of codes (cf. \cite{NRS01}, \cite{RaSl98}).
Usually	each new type of code is treated on an
individual basis.
We have recently found a
far-reaching generalization of Gleason's
result that applies simultaneously to arbitrary-genus
weight enumerators of self-dual codes over a very large class
of finite rings and modules.

The main result (Theorem \ref{mainII} below) can be
summarized as follows. Given a quasi-chain ring $R$ and
a notion of self-duality for codes over a left $R$-module $V$,
we construct a ``Clifford-Weil'' group $G$
such that the vector invariants of $G$ are spanned by the
full weight enumerators of self-dual isotropic codes,
and the polynomial invariants of $G$ are spanned by
the complete weight enumerators of these codes.

In the case of genus-$m$ weight enumerators (for $m \ge 1$) of
Type I binary self-dual codes, $G$ is the real
Clifford group ${\cal C} _m$ of our earlier paper \cite{NRS01}.
If $C$ is a Type II binary self-dual code, 
$G$ is the complex Clifford group ${\cal X}_m$
of \cite{NRS01}.
The case $m = 1$ gives the original Gleason theorem (except for
the specific identification of codes that generate the ring).
For self-dual codes over $\F_p$ containing the all-ones vector
(where $p$ is an odd prime) $G$ is the group
${\cal C}^{(p)}_m$ of \cite[Section 7]{NRS01}.

The proof of the main theorem will be the subject of our
forthcoming book \cite{BOOK}.
The proof is best carried out via a
categorical approach, using the concept of a form ring,
a generalization of the corresponding
concept from unitary $K$-theory \cite{HOM89}.
However, it is not necessary to understand this theory
to state and apply the theorem,
and so it seems worthwhile publishing a short version
that states the main result and gives some applications.
To illustrate the theorem we will
construct the Clifford-Weil groups for
generalized doubly-even codes over fields of characteristic 2
(Section \ref{Sec5}),
doubly-even codes over $\Z/2^f\Z$ (Section \ref{Sec6}),
and self-dual codes over
the noncommutative ring $\F_q + \F _q \, u$ where $u^ 2 = 0$ (Section \ref{Sec7}).
In Section \ref{Sec2} we define self-dual code,
isotropic code and the Type of a code.
Section \ref{Sec3} defines certain weight enumerators
and gives the main theorem, while Section \ref{Sec4}
studies the Clifford-Weil groups and their structure.
For further details the reader is referred to \cite{BOOK}.

\section{Self-dual isotropic codes.}\label{Sec2}

Throughout the paper,
$R$ will denote a ring (with unit element 1) and $V$ a left $R$-module.

\begin{defn}
A (linear) {\em code} $C$ of length $N$ over $V$ is
an $R$-submodule of the left $R$-module $V^N$.
In the classical language of coding theory, $V$
is the alphabet over which the code is defined, and
$R$ is the ground ring (so if $c \in C$ and $r \in R$,
$rc \in C$).
\end{defn}

Coding theory usually deals with codes over a {\em finite} alphabet $V$.
Therefore we will assume in the following that $V$
is a finite left $R$-module over the finite ring $R$.
Some parts of the theory carry over to infinite rings,
and so can be applied to self-dual lattices, for example
(see \cite{BOOK}). But
our explicit construction of the Clifford-Weil group
only applies to the finite case.

To express self-duality, we need a nonsingular bilinear form
$\beta $ on $V$. Since $V$ is finite, $\beta $ can be chosen
to have values in $\Q /\Z $. A bilinear form
$$\beta \in \Bil (V, \Q/\Z ) := \{ \beta : V\times V \to \Q/\Z \mid
\beta \mbox{ is $\Z $-bilinear }\} $$
is nonsingular if $v\mapsto \beta (v,\cdot ) $ is an isomorphism
of the abelian groups $V$ and $V^* = \Hom (V, \Q/\Z )$.
The following definition is then the natural generalization
of the usual notion of dual code (cf. \cite{MS77}).

\begin{defn}
Let $C \leq V^N$ be a code and $\beta \in \Bil (V, \Q/\Z )$ a
nonsingular bilinear form.
The {\em dual code} (with respect to $\beta $) is
$$C^{\perp} := \{ x=(x_1,\ldots, x_N) \in V^N \mid \sum _{i=1}^N \beta (x_i,c_i ) = 0 
\mbox{~for~all~} c=(c_1,\ldots, c_N ) \in C \} ~.$$
$C$ is called {\em self-dual} (with respect to $\beta $) if
$C=C^{\perp }$.
\end{defn}

To express certain additional constraints
on the code (that weights are divisible by 4
for binary codes, or that the code contains the
all-ones vector ${\bf 1} = (1,\ldots, 1)$, etc.)
we use quadratic mappings, which we define to be
sums of quadratic forms and linear forms on $V$, with values in $\Q / \Z$:

\begin{defn}
Let
$\Quad _0(V, \Q/\Z) :=  \{ \phi : V\to \Q/\Z \mid \phi(0) = 0$ and
$$\phi(x+y+z) - \phi(x+y) - \phi(x+z) - \phi(y+z) + \phi(x) + \phi(y) + \phi(z) = 0 \} ~.$$
Let $\Phi \subset \Quad _0 (V,\Q/\Z) $ and let $C \leq V^N$ be a code.
Then $C$ is called {\em isotropic} (with respect to $\Phi $)
if $\phi (c) := \sum _{i=1}^N \phi (c_i) = 0$ for all $c\in C$ and $\phi \in \Phi $.
\end{defn}

\begin{remarks}
Let $V$ be a left $R$-module. Then

(a)
$\Bil (V, \Q/\Z ) $ is a right $(R\otimes R)$-module via
$\beta (r\otimes s) (v,w) = \beta (rv,sw) $ for all
$\beta \in \Bil (V,\Q/\Z )$, $v,w \in V$, $r,s\in R$.

(b)
$\Bil (V,\Q/\Z) $ has a natural involution
$\tau $ defined by $(\beta ^\tau)(v,w) := \beta (w,v) $ for all
$\beta \in \Bil (V,\Q/\Z )$, $v,w \in V$.

(c)
For $r \in R$ define $[r] \in \End(\Quad _0(V,\Q/\Z ))$ by
$(\phi [r])(v) = \phi (rv) $ for $\phi \in \Quad _0(V,\Q/\Z )$,
$v \in V$.
This ``action'' satisfies
$ [rs] = [r][s] $ and
$$[r+s+t] -[r+s]-[r+t]-[s+t] +[r]+[s]+[t] = 0~,$$ for all $r,s,t \in R$.
We describe this situation by saying that we have made
the group $\Quad _0(V,\Q/\Z )$ into an {\em $R$-qmodule}.

(d)
The mapping $[-1]$ is an involution on $\Quad _0(V,\Q/\Z )$.

(e)
There is a mapping $\leftBra \rightBra : \Bil (V,\Q/\Z) \to \Quad _0 (V,\Q/\Z)$
defined by $\leftBra \beta \rightBra (v) := \beta (v,v) $ for all $\beta \in
\Bil (V,\Q/\Z )$, $v\in V$.

(f)
There is a mapping $\lambda : \Quad _0 (V,\Q/\Z) \to \Bil (V,\Q/\Z)$
defined by $(\lambda (\phi )) (v,w) := \phi (v+w) - \phi (v) - \phi (w) $ for all $\phi \in \Quad _0 (V,\Q/\Z )$, $v,w \in V$.

(g)
Both mappings $\lambda $ and $\leftBra \rightBra $ are $R$-qmodule homomorphisms,
where $\Bil (V, \Q/\Z ) $ is regarded as an $R$-qmodule
via $\beta [r] (v,w) := \beta (rv,rw)$ for all
$\beta \in \Bil (V,\Q/\Z )$, $v,w \in V$, $r\in R$.

(h)
For all $\beta \in \Bil (V, \Q/\Z )$ and $\phi \in \Quad _0 (V,\Q/\Z )$ we have
\begin{align}
\leftBra \beta ^\tau\rightBra &= \leftBra \beta \rightBra ~, \nonumber \\
\lambda(\phi)^\tau &= \lambda(\phi) ~, \nonumber \\
\lambda(\leftBra \beta \rightBra ) &= \beta + \beta ^\tau ~, \nonumber \\
\phi[r+s] - \phi [r] - \phi [s] &= \leftBra \lambda(\phi)(r\otimes s)\rightBra ~. \nonumber
\end{align}

(i)
For all $\beta \in \Bil (V, \Q/\Z )$ and $\phi \in \Quad _0 (V,\Q/\Z )$ we have
$$\lambda ( \leftBra \lambda (\phi ) \rightBra ) = 2 \lambda (\phi )
\mbox{ and }
\leftBra \lambda (\leftBra \beta \rightBra ) \rightBra = 2 \leftBra \beta \rightBra ~.$$
\end{remarks}

We call $\beta \in \Bil (V,\Q/\Z )$ {\em admissible}
if $\beta $ is nonsingular and the $(1\otimes R)$-submodule
of $\Bil (V,\Q/\Z )$ given by
$$ M:= \beta (1\otimes R) = \{ \beta _r \mid r\in R \}, 
\mbox{~where~}\beta (1\otimes r) = \beta_r \mbox{~is defined by~} 
\beta _r (v,w) = \beta ( v,rw ) ~,$$
is closed under $\tau $ and isomorphic to
$R$, i.e. if $$\psi: r \mapsto \beta _r$$ defines an
isomorphism of right $R$-modules
$\psi : R_R \to M_{1\otimes R}$.

Note that any admissible $\beta $ defines an anti-automorphism
$J$ of $R$ by $r\mapsto r^J$, where
$$\beta (rv,w) = \beta (v, r^J w ) \mbox{~for~all~} v,w \in V .$$
Let $\epsilon \in R $ be defined by
$\beta (v,w) = \beta (w, \epsilon v) $ for all $v,w\in V$.
Then $\epsilon ^J r^{(J^2)} \epsilon = r$ for all $r\in R$.
In particular $\epsilon ^J \epsilon = 1$ and since $R$ is finite,
$\epsilon $ is a unit.

\begin{defn}
The quadruple $\rho := (R,V,\beta, \Phi )$ is called a
({\em finite representation of a}) {\em form ring}
if $R$ is a finite ring, $V$ is a finite left $R$-module,
$\beta \in \Bil (V,\Q/\Z )$ is an admissible bilinear form,
and $\Phi \leq \Quad _0 (V,\Q/\Z )$ is a sub-$R$-qmodule such that
$\leftBra M \rightBra \leq \Phi $ and
$\lambda (\Phi ) \leq M $, where $M:= \beta (1\otimes R)$.
Then if $C\leq V^N $ for some integer $N \ge 1$
is self-dual with respect to $\beta $ and isotropic
with respect to $\Phi $, $C$ is called a (self-dual isotropic)
{\em code of Type $\rho $}.
\end{defn}

\begin{remarks}
(a)
The definition of Type given in \cite{BOOK} 
uses a more abstract definition of a form ring.
For our purposes here it is enough to work with the
above concrete realization.

(b)
The definition of Type given above (and the still more general
version given in \cite{BOOK}), are far-reaching generalizations of the
notion of Type used in \cite{MS77}, \cite{RaSl98}.

(c)
The usual notion of form ring
\cite{HOM89} corresponds directly to the case when $\lambda$ is injective;
it is also closely connected to the dual case when 
$\leftBra\rightBra$ is surjective.

\end{remarks}

\section{Weight enumerators and Clifford-Weil groups.}\label{Sec3}

There are many kinds of weight enumerators of codes
(see e.g. \cite{RaSl98}).
We introduce here only the complete weight enumerator, although analogues of
our main theorem hold for full weight enumerators
and other appropriate symmetrizations of the full weight enumerator.

\begin{defn}
Let $C\leq V^N$ be a code.
The
{\em complete weight enumerator}
of $C$ is
$\cwe(C) := \sum _{c\in C}
\prod _{i=1}^N x_{c_i}
 \in \C [x_v \mid v\in V ] $.
\\
The {\em genus-$m$ complete weight enumerator} of $C$ is
$$\cwe_m(C) := \sum _{(c^{(1)},\ldots, c^{(m)})\in C^m}
\prod _{i=1}^N x_{(c^{(1)}_i,\ldots,c^{(m)}_i) }
 \in \C [x_v \mid v\in V ^m] ~.$$
\end{defn}

\begin{rems}{\label{cwem}}
The genus-$m$ complete weight enumerator of $C$ can be obtained from the complete
weight enumerator of $C(m):=C\otimes R^m \leq V^N\otimes R^m \cong V^{mN}$.
For if we identify $V^{mN}$ with $(V^m)^N$ and consider
$C(m)$ as a code in $(V^m)^N$, then
$\cwe_m(C) = \cwe(C(m))$.
\end{rems}

Note that $C(m)$ is a code for which the ground ring is $\Mat _m(R)$, the ring of $m\times m$-matrices with entries in $R$.
This is one of the main reasons why we allow 
noncommutative ground rings:
even when we consider
genus-$m$ weight enumerators for classical binary
codes, $\Mat _m(\F_2) $ arises naturally as the ground ring.

Let $\rho = ( R,V,\beta, \Phi )$ be a form ring.
Let $C\leq V^N$ be a self-dual isotropic code of Type $\rho $.
Then the complete weight enumerator
$\cwe (C)$ is invariant under the substitutions
$$\rho (r) : x_v \mapsto x_{rv} \ , \ \mbox{~for~all~} r\in R^*
\hspace{1cm} \hfill{\mbox{(since $C$ is a code)}} ~, $$
$$\rho (\phi ) : x_v \mapsto \exp(2\pi i \phi (v) ) x_{v} \ ,
 \mbox{~for~all~} \phi \in \Phi
\hspace{1cm} \hfill{\mbox{(since $C$ is isotropic)}} ~,$$ as well as the
MacWilliams transformation (cf. \cite{MS77}, \cite{RaSl98}): \\
$$ h : x_v \mapsto \sqrt{|V|}^{-1} \sum _{w\in V } \exp (2\pi i \beta (w,v ) ) x_{w} \hspace{1cm} \hfill{\mbox{(since $C=C^{\perp}$).}} $$

Subsidiary transformations can be derived from $h$.
If $R=\Z/6\Z$, for example, we also obtain the MacWilliams transformations
modulo 2 and modulo 3.
For general rings $R$, one can similarly construct further
MacWilliams transformations using symmetric idempotents,
as we now demonstrate.

\begin{example}{\label{symmetricidempotent}}
Let $V'=\iota R$ for some idempotent $\iota \in R$. Then $V'$ admits a
nonsingular $J$-Hermitian form if and only if there is an
isomorphism of right $R$-modules
$$\kappa : \iota R\cong \iota^J R ~, $$
in which case we say that the idempotent is {\em symmetric}.
Note that any such isomorphism $\kappa $
has the form
$$
\kappa(\iota x)=v_\iota x,\ \kappa^{-1}(\iota^J x)=u_\iota x,
$$
where $u_\iota\in \iota R \iota^J$ and $v_\iota\in \iota^J R \iota$
satisfy
$
u_\iota v_\iota = \iota,\ v_\iota u_\iota = \iota^J.
$
\end{example}

If $\iota = u_{\iota } v_{\iota } $ is a symmetric idempotent in $R$, then
$$h _{\iota,v_{\iota}} : x_v \mapsto \sqrt{|\iota V|}^{-1} \sum _{w\in \iota V } \exp (2\pi i \beta (w, v_{\iota } v ) ) x_{w+(1-\iota) v} $$
is the partial MacWilliams transformation corresponding to $\iota V$.

\begin{defn}
Let $\rho = (R,V,\beta, \Phi )$ be a form ring.
Then the associated {\em Clifford-Weil group} is
$${\cal C}(\rho ) := \langle \rho (r), \rho(\phi ), h_{\iota,v_{\iota} } \mid
r\in R^*, \phi \in \Phi, \iota \mbox{ symmetric idempotent in } R \rangle ~,
$$
a subgroup of $\GL _{|V|} (\C )$.
\end{defn}

The above discussion shows
that the complete weight enumerator of a self-dual
isotropic code of Type $\rho $ is invariant under the action of
the group ${\cal C}(\rho )$.
If $G$ is any subgroup of $\GL _{n} (\C )$ we let
$$\Inv (G) := \{ p \in \C [x_1,\ldots, x_n ] \mid p(gX) = p(X) \mbox{~for~all~}
g\in G \} $$
denote the invariant ring of $G$.

Now we can state our main theorem, part (i) of which is clear from
the above considerations and the MacWilliams identities.
We cannot at present prove part (ii) for arbitrary finite rings, but
need to make some additional assumptions, for example that
$R$ is a {\em chain ring}
(i.e. the left ideals in $R$ are linearly ordered by inclusion)
or, more generally, a {\em quasi-chain ring},
by which we mean a direct product of matrix rings over chain rings.

\begin{thm}\label{mainII}
Let $\rho = (R,V,\beta, \Phi )$ be a form ring.
\begin{itemize}
\item[(i)] If $C\leq V^N$ is a self-dual isotropic code
then $\cwe (C) \in \Inv ({\cal C}(\rho )) $.
\item[(ii)] If $R$ is a finite quasi-chain ring
then $\Inv ({\cal C}(\rho )) $
is spanned by complete weight enumerators of self-dual isotropic codes of
Type $\rho $:
$$\Inv ({\cal C}(\rho )) = \langle \cwe (C) \mid C \mbox{ self-dual, isotropic
code in } V^N, N \ge 1 \rangle .$$
\end{itemize}
\end{thm}

To deal with higher-genus weight enumerators
we introduce the {\em associated Clifford-Weil group}
${\cal C}_m(\rho ):= {\cal C}(\rho \otimes R^m)$
{\em of genus $m$}. By Remark \ref{cwem}
the invariant ring of ${\cal C}_m(\rho )$
is spanned by the genus-$m$ weight enumerators of self-dual codes
of Type $\rho $.

Let $\rho := (R,V,\beta,\Phi)$ be a form ring.
By Morita theory (see \cite{BOOK}),
this corresponds to a unique form ring
$$\Mat _m(\rho) = \rho \otimes R^m := (\Mat_m(R),V\otimes R^m,\beta ^{(m)},\Phi _m) ~,$$
which we call a {\em matrix ring for the form ring $\rho $}.
Here $\beta ^{(m)}$ is the bilinear form
on $V^m = V\otimes R^m$ (admissible for $\Mat_m(R)$) defined by
$$\beta ^{(m)} ((v_1,\ldots, v_m), (w_1,\ldots, w_m) ) =
\sum _{i=1}^m \beta (v_i,w_i ) ~.$$
$\beta ^{(m)}$ generates the $\Mat _m(R)$-module
$M_m := \beta ^{(m)} (1 \otimes \Mat _m(R) )$, and
$$\Phi _m = \left\{ \left( \begin{array}{cccc}
\phi _1 & m_{12} & \ldots & m_{1m} \\
 & \ddots & \ddots & \vdots \\
& & \ddots & m_{m-1,m} \\
& & & \phi _m \end{array} \right) \mid \phi_1,\ldots, \phi _m \in \Phi, m_{ij} \in M \right\}  $$
is a set of upper triangular matrices,
where
$$\left( \begin{array}{cccc}
\phi _1 & m_{12} & \ldots & m_{1m} \\
 & \ddots & \ddots & \vdots \\
& & \ddots & m_{m-1,m} \\
& & & \phi _m \end{array} \right) (v_1,\ldots, v_m)
= \sum _{i=1}^m \phi _i (v_i) + \sum _{i<j} m_{ij} (v_i,v_j) ~.$$
One easily sees that $\Phi _m$ is a sub-$\Mat _m(R)$-qmodule of
$\Quad _0(V\otimes R^m)$.
The involution ${J_m}$ on $\Mat _m(R)$
acts as componentwise application of $J$ followed by
transposition.
The unit $\epsilon _m \in \Mat _m(R)$ is the scalar matrix
$\epsilon I_m$.

\section{The structure of the Clifford-Weil groups.}\label{Sec4}

The group ${\cal C}(\rho )$ is a projective representation of a
so-called {\em hyperbolic co-unitary group} ${\cal U}(R,\Phi )$,
which we will define in terms of $R$, the involution
$J$ and the $R$-qmodule $\Phi $.
${\cal U}(R,\Phi )$ is an extension of the linear
$R$-module $\ker(\lambda ) \oplus \ker (\lambda )$ by a
certain group $G$ of $2\times 2$ matrices over the ring $R$.
This is one of two possible
extensions of the unitary $K$-theoretic notion of ``hyperbolic unitary
group'' \cite{HOM89} that applies to our generalized definition of form ring.

Let $\rho :=(R,V,\beta, \Phi )$ be a form ring and
let $M:=\beta (1\otimes R)$.
Then, using the above construction of the $2\times 2$-matrix ring for a form ring,
$\Phi _2$ is a $\Mat _2(R)$ sub-qmodule of $\Quad _0(V\otimes R^2)$.
%$$\Phi _2 = \{ \left( \begin{array}{cc} \phi _1 & m \\ & \phi _2
%\end{array} \right) \mid \phi _1, \phi _2 \in \Phi, m\in M \} .$$
%Imitating the formal matrix multiplication
%$A^{J_2} \phi A $ we get the action of $\Mat_2(R)$ on $\Phi _2$ by
%$$
%\left( \begin{array}{cc} \phi _1 & m \\ & \phi _2 \end{array} \right)
%[ \left( \begin{array}{cc} a & b \\ c & d \end{array} \right) ]
%:=
%\left( \begin{array}{cc} \phi' _1 & m' \\ & \phi' _2 \end{array} \right)
%$$
%where
%\begin{align*}
%\phi' _1 &:= \phi _1 [a] + \phi _2 [c] + \leftBra m (a\otimes c) \rightBra
%~, \nonumber \\
%m' &:= \lambda(\phi _1) (a\otimes b) + m (a\otimes d) + \lambda(\phi _2 ) (c\otimes d) + m ^\tau (c\otimes b) ~, \nonumber \\
%\phi' _2 &:= \phi _1 [b] + \phi _2 [d] + \leftBra m (b\otimes d) \rightBra ~. \nonumber
%\end{align*}
%One can check that this turns $\Phi _2$ into a $\Mat_2(R)$-qmodule.
In particular, $\Phi _2$ is a module for the unit group $\GL _2 (R) $ of $\Mat_2(R)$,
and we can form the semi-direct product $\GL _2(R) \ltimes \Phi _2$
of which ${\cal U}(R,\Phi )$ will be a subgroup.
Applying $\lambda $ and $\psi $ componentwise, we get mappings
$\lambda _2 : \Phi _2 \to \Mat _2(M) $ defined by
$\lambda _2(\left( \begin{array}{cc} \phi _1 & m \\ & \phi_2 \end{array}\right)) :=
\left( \begin{array}{cc} \lambda(\phi _1) & m \\ m ^\tau & \lambda(\phi _2) \end{array} \right)$
and
$\psi _2 : \Mat _2 (R) \to \Mat_2(M) $.

Then
${\cal U}(R,\Phi ) = $ $$ \{
(\left( \begin{array}{cc} a & b \\ c & d \end{array} \right),
\left( \begin{array}{cc} \phi _1 & m \\ & \phi_2 \end{array}\right))
\in \GL _2(R) \ltimes \Phi _2
\mid
\psi _2 (
\left( \begin{array}{cc} c^Ja & c^Jb \\ d^Ja-1 & d^Jb \end{array} \right) )
%\left( \begin{array}{cc} \psi(c^Ja) & \psi(c^Jb) \\ \psi(d^Ja-1) & \psi(d^Jb) \end{array} \right)
=
\lambda _2 (\left( \begin{array}{cc} \phi _1 & m \\ m ^\tau & \phi _2 \end{array} \right)) \} .$$

\begin{remark}{\label{isounihyp}}
To describe the isomorphism type of ${\cal U}(R,\Phi )$,
note that the projection
$$\pi: {\cal U}(R,\Phi ) \rightarrow \GL_2(R), \ (u,\phi ) \mapsto u ~,$$
defines a group homomorphism. The kernel of $\pi $ is the set of
all $(1,\phi )\in {\cal U}(R,\Phi )$, i.e.
$$\ker (\pi ) = \{ (1,
\begin{pmatrix} \phi _1 & 0 \\ & \phi _2 \end{pmatrix} ) \mid
\phi _1,\phi _2 \in \Phi,\ \lambda (\phi _1) = \lambda (\phi _2) =0 \} ~,$$
which is naturally isomorphic to $\ker(\lambda ) \times \ker(\lambda )$.
The image of $\pi $ is
$$\{ \begin{pmatrix} a & b \\ c & d \end{pmatrix} \in \Mat_2(R) \mid
\left( \begin{array}{cc} a^J & c^J \\ b^J & d^J
\end{array} \right)
\begin{pmatrix} 0 & 0 \\ 1 & 0 \end{pmatrix}
\left( \begin{array}{cc} a & b \\ c & d
\end{array} \right) -
\begin{pmatrix} 0 & 0 \\ 1 & 0 \end{pmatrix} \in \psi_2^{-1}(\lambda _2(\Phi _2)) \} ~.$$
In many important examples it is easy to describe the image of $\psi _2 ^{-1} \circ \lambda _2$
(for example, this may consist of
all symmetric, skew-symmetric or hermitian elements)
and so leads to an isomorphism of
${\cal U}(R,\Phi ) / \ker(\pi) $ with (a subgroup of) a classical group.
\end{remark}

\begin{example}
Let $R$ be one of the finite simple rings with involution
shown in the following table,
where $\Mat _n(\F _q)$ denotes the ring of $n \times n$ matrices over
the field $\F_q$, with transposition denoted by ${\ttt}$.
Then for an appropriate center $Z$,
$${\cal C}(\rho ) = Z \,.\, {\cal U}(R,\Phi ) = Z \,.\, (\ker(\lambda ) \oplus \ker(\lambda )) \,.\, {\cal G}(R,\Phi ) ~,$$
where ${\cal G}(R,\Phi )$ is the classical group shown in the last column:
\begin{center}
\begin{tabular}{|c|c|c|c|}
\hline
$R$ & $J$ & $\epsilon $ & ${\cal G}(R,\Phi )$ \\
\hline
& & & \\
$\Mat _n(\F _q) \oplus \Mat _n(\F _q) $ & $(r,s)^J = (s^{\ttt},r^{\ttt}) $ & 1 & $\GL _{2n}(\F_q )$ \\
\hline
& & & \\
$\Mat _n(\F _{q^2}) $ & $r^J = (r^q)^{\ttt } $ & 1 & $U _{2n}(\F_{q^2} )$ \\
\hline
& & & \\
$\Mat _n(\F _q), \ q=p^m, \ p>2 $ & $r^J = r^{\ttt } $ & 1 & $\Sp _{2n}(\F_{q} )$ \\
\hline
& & & \\
$\Mat _n(\F _q), \ q=p^m, \ p>2 $ & $r^J = r^{\ttt } $ & $-1$ & $O^+ _{2n}(\F_{q} )$ \\
\hline
& \multicolumn{2}{l|}{} & \\
$\Mat _n(\F _q), \ q=p^m, \ p=2 $ & \multicolumn{2}{l|}{$\psi^{-1}(\lambda (\Phi)) = \{ r\in R \mid r^J = r \}$} & $\Sp_{2n}(\F_{q} )$ \\
\hline
& \multicolumn{2}{l|}{} & \\
$\Mat _n(\F _q), \ q=p^m, \ p=2 $ & \multicolumn{2}{l|}{$\psi^{-1}(\lambda (\Phi)) = \{ 0 \}$} & $O^+ _{2n}(\F_{q} )$ \\
\hline
\end{tabular}
\end{center}
\end{example}

Recall that a ring $R$ is semiperfect
(cf.  \cite[page 346]{Lam01})
if $R/\rad R$ is semisimple
and idempotents of $R/\rad R$ lift to idempotents of $R$;
in particular, all finite rings are semiperfect.
One can show the following:

\begin{thm} (See \cite{BOOK}.)
If $R$ is semiperfect, then the hyperbolic co-unitary group ${\cal U}(R,\Phi )$ is 
generated by the following elements:
\[
d(u,\phi) := (\begin{pmatrix} u^{-J}&u^{-J} \psi ^{-1} (\lambda(\phi))\\0&u\end{pmatrix},
\begin{pmatrix} 0&0\\&\phi\end{pmatrix})
\]
for all $u\in R^*$, $\phi \in \Phi $, and
\[
H_{\iota,u_\iota,v_\iota}
:=
(\begin{pmatrix}
1-\iota^J&v_\iota\\
-\epsilon^{-1} u_\iota^J & 1-\iota\end{pmatrix}
,
\begin{pmatrix}
0&\psi(-\epsilon\iota )\\
& 0
\end{pmatrix}
)
\]
for all symmetric idempotents $\iota \in R$.
\end{thm}

Then the projective representation ${\cal U}(R,\Phi )\to {\cal C}(\rho )$ is
defined by
$$d(u,\phi ) \mapsto \rho (u) \rho (\phi )$$ and
$$H_{\iota,u_\iota,v_\iota} \mapsto h_{\iota,v_\iota} .$$

Hence ${\cal C}(\rho )$ is an extension
$$ 1 \to Z \to {\cal C}(\rho ) \to {\cal U}(R,\Phi ) \to 1 $$
where $Z$ consists of scalar matrices (since the projective
representation can be seen to be irreducible).
If Theorem \ref{mainII} holds for $\rho $, then by invariant theory
$Z\cong Z_l$ is cyclic of order $l$ where
$l = \gcd \{ N \mid \mbox{~there exists~} C\leq V^N \mbox{ of Type } \rho \} $.
In fact, for arbitrary $\rho$ (not necessarily satisfying
Theorem  \ref{mainII}), one can show that
\begin{align*}
|Z| & ~=~  \gcd \{N \ge 1 \mid \mbox{there~exists~} C\le V^N \mbox{~of~Type~} \rho\} \\
 & ~=~  \min \{N \ge 1 \mid \mbox{there~exists~} C\le V^N \mbox{~of~Type~} \rho\} ~.
\end{align*}

\section{Doubly-even euclidean self-dual codes over $\F_{2^f}$.}\label{Sec5}
This and the next two sections will illustrate the above theory.
For further details about this first section see \cite{NQRS}.
In \cite{Queb91} Quebbemann defines the notion of an
even code over the field $k:=\F _{2^f}$ as follows.
A code $C\leq k ^N $ is called {\em even} if
\beql{Eq1}
\sum _{i=1}^N c_i = 0 \mbox{ and } \sum _{i<j} c_i c_j =0, \ \
\mbox{~for~all~} c\in C ~.
\eeq
It is easy to see that even codes are self-orthogonal
with respect to the usual bilinear form $\sum _{i=1}^N c_i c'_i $.
Moreover, identifying $k$ with $\F_2^f$ using a
self-complementary (or trace-orthonormal) basis, even codes
remain even over $\F _2$.
A self-dual even code in this sense is called a
{\em generalized doubly-even self-dual code}.
If $f=1$, even codes are precisely the classical doubly-even
(or Type II) binary codes.
If $f=2$, they are the Type II codes over $\F_4$ considered in \cite{GPSA02}.

The Type of these codes can be specified in the language of
form rings as follows.
Let $R=\F _{2^f}$, $V=R$ and $\beta : V\times V \to \frac{1}{2} \Z/\Z $ be
defined by $\beta (x,y) := \frac{1}{2} \tr (xy) $, where
$\tr $ denotes the trace from $\F _{2^f} $ to $\F_2 \cong \Z/2\Z $.
Then $\beta $ is admissible and $M:=\beta (1\otimes R) = \{ \beta _a := \beta (1\otimes a) \mid a \in R \} $.

The quadratic forms will take values modulo $4$.
Let $O:= \Z_2[\zeta _{2^f-1}]$ be the ring of integers in the unramified
extension of degree $f$ of the 2-adic numbers.
Then $R \cong O/2O$. If $x\in R$, $x^2$ is uniquely determined modulo 4,
so squares of elements of $R$ can be considered as elements of $O/4O$.
The usual trace $\Tr : O \to \Z $ maps $4O $ into $4\Z$.
For $a\in R$ we define
$$\phi _a : V \to \frac{1}{4} \Z/ \Z, \ 
\phi _a(x):= \frac{1}{4} \Tr (a^2x^2) \in \Quad _0(V,\Q/\Z) ~,$$
and let $\Phi := \{ \phi _a \mid a\in R \}$.
Then $(R,V,\beta, \Phi )$ is a form ring.

\begin{thm}
Codes of Type $\rho $ are exactly 
the generalized doubly-even self-dual
codes in the sense of \eqn{Eq1}.
\end{thm}

\begin{proof}
Let $C\leq \F_{2^f}^N$ be an even code in the sense of \eqn{Eq1}.
Since $\lambda $ is surjective, it is enough to show that
$\sum _{i=1}^N \phi _a(c_i) = 0 $ for all $c\in C$.
Now $\sum _{i=1}^N c_i = 0$,
therefore as an element of $O/4O$ the square
\beql{Eq2}
(\sum _{i=1}^N c_i)^2 =\sum _{i=1}^N c_i^2 + 2\sum _{i<j} c_ic_j = 0 \in O/4O ~.
\eeq
Since $\sum _{i<j} c_i c_j = 0$ it follows that
$\sum _{i=1}^N c_i^2 = 0 \in O/4O$.

To obtain the other inclusion, let $C$ be a code of Type $\rho $.
By the nondegeneracy of the trace form, $\sum _{i=1}^N c_i^2 = 0 $ in
$O/4O$.
Therefore by \eqn{Eq2},
$ (\sum _{i=1}^N c_i)^2 \equiv 0 \pmod{2 O} $ and hence also
$\sum _{i=1}^N c_i = 0 \in \F _{2^f}$.
Then
$ (\sum _{i=1}^N c_i)^2 \equiv 0 \pmod{4 O} $ and \eqn{Eq2} implies that
 $\sum _{i<j} c_i c_j = 0$.
\end{proof}

We next compute the
Clifford-Weil groups of arbitrary genus for these codes.

\begin{thm}
Let ${\cal C}_m(\rho )$ be the Clifford-Weil group of genus $m$
corresponding to
the form ring $\rho $ above.
Then $${\cal C}(\rho ) \cong Z\,.\,(k^m\oplus k^m)\,.\,\Sp_{2m}(k)
\cong Z \,.\, 2^{2mf} \,.\, \Sp _{2m}(2^f) ~,$$
where $Z \cong Z_4 $ if $f:=[k:\F_2 ] $ is even, 
$Z \cong Z_8 $ if $f:=[k:\F_2 ] $ is odd.
\end{thm}

\begin{proof} From Section \ref{Sec4},
${\cal C}_m(\rho ) $ has an epimorphic image
$${\cal U}(\Mat _m(R),\Phi _m ) \cong (k^m\oplus k^m) . \Sp_{2m}(k) ~.$$
The kernel $Z$ of this epimorphism is a cyclic group
consisting of scalar matrices.
Since the invariant ring of ${\cal C}_m(\rho )$ is
spanned by weight enumerators of self-dual isotropic codes $C$,
the order of $Z$ is the greatest common divisor of the
lengths of these self-dual isotropic codes.
Since the codes are self-dual, they all contain the
all-ones vector ${\bf 1}$. This vector spans an isotropic
$k$-space if and only if the length $N$ of the code is divisible by 4.
Therefore $|Z|$ is divisible by $4$.

First consider the case when $f$ is even.
The code $Q_4$ with generator matrix
$$
\begin{bmatrix} 1 & 1 & 1 & 1 \\ 0 & 1 & \omega & \omega ^2 \end{bmatrix} ~,
$$
where $\omega \in \F_4 \setminus \F_2 $, is an isotropic self-dual
code of length 4 over $\F_4$.
Extending scalars, one gets an isotropic self-dual code
$k\otimes _{\F _4} Q_4 $ of length 4 for all $k$ which are
of even degree over $\F_2$.
Hence in this case $Z \cong Z_4 $.

If $f$ is odd, then any self-dual isotropic code
$C\leq k^N$ yields a doubly-even self-dual binary code
$\tilde{C} \leq \F _2^{fN}$. The length $fN$ of $\tilde{C}$ is
necessarily divisible by 8. Since $f$ is odd, this implies that $N$
is divisible by 8 and hence in this case $Z \cong Z_8 $.
\end{proof}

{\bf The case $k=\F_2$.} We obtain
the classical doubly-even binary codes;
the higher-genus Clifford-Weil groups are the complex Clifford groups of
\cite{NRS01}.

{\bf The case $k=\F_4$.}
Let $k=\F_4 = \{ 0,1,\omega, \omega ^2 \}$.
Then
$$G:={\cal C}(\rho ) \cong (Z_4  \zentr  D_8  \zentr  D_8 ) \,.\, \mbox{Alt}_5 $$
(where $ \zentr $ denotes a central product) is generated by
$$\rho (\omega ) =
 \left( \begin{array}{cccc}
1 & 0 & 0 & 0 \\
0 & 0 & 0 & 1 \\
0 & 1 & 0 & 0 \\
0 & 0 & 1 & 0 \end{array} \right) ,
~h=\frac{1}{2}
 \left( \begin{array}{rrrr}
1 & 1 & 1 & 1 \\
1 & 1 & -1 & -1 \\
1 & -1 & -1 & 1 \\
1 & -1 & 1 & -1 \end{array} \right) ,
$$
and
$\rho (\phi ) = \diag (1,-1,i,i ) $.
$G$
is a subgroup of index 2 in the
complex reflection group $G_{29}$ 
(No. 29 in \cite{ShTo54}).
The Molien series of $G$  (cf. \cite{MS77}, \cite{RaSl98}) is
$$\frac{1+t^{40}}{(1-t^4)(1-t^8)(1-t^{12})(1-t^{20})} ~.$$
Primary invariants of $G$ (which generate the invariant ring of
$G_{29}$) can be taken to be the complete weight enumerators of the extended
quadratic residue codes
$Q_N$ of lengths $N=4,8,12,20$ over $\F_4$.
%$Q_i:=\tilde{QR}(\F_4,i-1)$ of length $i=4,8,12,20$.
The elements in $G_{29} \setminus G$ act as the Frobenius
automorphism $\omega \mapsto \overline{\omega }$
on the weight enumerators of these codes.
A weight enumerator $p_C$ is invariant under $G_{29}$
if $p_C = \overline{p_C}$.
To get the full invariant ring of $G$, a further code is needed,
a self-dual even code $C_{40}$ of length 40 over $\F_4 $,
for which the complete weight enumerator is not invariant
under the Frobenius automorphism.
Such a code is constructed in \cite{BeCh02}.

{\bf The case $k=\F_8$.}
In this case we will just give the Molien series
of ${\cal C}(\rho )$. This can be written as
$f(t)/g(t)$, where
\begin{eqnarray*}
& f(t) := & 1+ 5t^{16} + 77t^{24} + 300t^{32} + 908t^{40} + 2139t^{48} + 3808t^{56} + 5864t^{64} \\
&&{} ~ + 8257t^{72} + 10456t^{80} + 12504t^{88} + 14294t^{96} + 15115t^{104} \\
&&{} ~ + 15115t^{112} + 14294t^{120} + 12504t^{128} + 10456t^{136} + 8257t^{144} \\
&&{} ~ + 5864t^{152} + 3808t^{160} + 2139t^{168} + 908t^{176} + 300t^{184} \\
&&{} ~ + 77t^{192} + 5t^{200} + t^{216}
\end{eqnarray*}

and
$$ g(t) := (1-t^8)^2(1-t^{16})^2(1-t^{24})^2(1-t^{56})(1-t^{72}) ~.$$
It would clearly be hopeless to attempt to find 
codes whose weight enumerators generate this ring.
This phenomenon is typical (compare \cite{HuSl79}):
it is the exception rather than the rule
for these rings of invariants to have a simple
description in terms of codes.

\section{Doubly-even self-dual codes over $\Z/2^f\Z $.}\label{Sec6}

Let $R:=\Z/2^f\Z$ and let $C \leq R^N$ be a code of length $N$.
Then the dual code
$C^{\perp } := \{ x\in R^N \mid \sum _{i=1}^N x_ic_i = 0, 
\mbox{~for~all~} c\in C \}$.
$C$ is called {\em doubly-even} if
$\sum _{i=1}^N c_i^2 \equiv 0 \pmod{2^{f+1}} $.
To describe the class of doubly-even self-dual codes over $R$
in the language of form rings let $V:= R$,
and define $\beta : V\times V \to \frac{1}{2^f}\Z /\Z $ by $\beta (x,y) := \frac{1}{2^f} xy $
and $\phi _0: V \to \frac{1}{2^{f+1}} \Z / \Z $ by
$\phi_0 (x) := \frac{1}{2^{f+1}} x^2$.
Let $\Phi $ be the $R$-qmodule generated by $\phi _0 $ and 
define the form ring
$$\rho _a:= (R,V,\beta, \Phi ).$$

To express the additional property that $C$ contains the all-ones vector,
let $\varphi  : V \to \frac{1}{2^{f}} \Z / \Z $ be defined by
$\varphi (x) := \frac{1}{2^{f}} x$ and let $\Phi _0 $ be the $R$-qmodule spanned
by $\phi _0 $ and $\varphi $.
Then define the form ring
$$\rho _b:= (R,V,\beta, \Phi _0 ).$$
Note that $\ker (\lambda _a) = \langle 2^f \phi _0 \rangle
\cong \Z/2\Z $, whereas
$\ker (\lambda _b) =\langle \varphi  \rangle \cong R $, since $2^{f-1} \varphi  = 2^f \phi _0 $.
Since the involution on $R$ is trivial and $\epsilon = 1$, we find that
$$\pi ({\cal U}_m(R,\Phi ) )= \{ A \in \Mat_{2m}( R) \mid
A^{tr} \begin{pmatrix} 0 & -I_m \\ I_m & 0 \end{pmatrix} A =
\begin{pmatrix} 0 & -I_m \\ I_m & 0 \end{pmatrix} \} \cong \Sp_{2m}(R) .$$

\begin{thm}
For $x\in \{a,b\}, $ let
$l_x = \gcd \{ N \mid \mbox{~there exists~a~code~} C\leq V^N \mbox{ of Type } \rho _x \}$.
The Clifford-Weil groups of genus $m$ are extensions
$${\cal C}(\rho _a) \cong Z_{l_a}\,.\, (Z_2^m \times Z_2^m)\,.\, \Sp_{2m}(R) $$
and
$${\cal C}(\rho _b) \cong Z_{l_b}\,.\, (R^m \times R^m)\,.\, \Sp_{2m}(R) $$
\end{thm}

As we show in \cite{BOOK}, $l_a = 8$ and $l_b = \max \{2^{f+1}, 8\} $.

{\bf The case $f=2$.}
We identify $R$ with $\Z/4\Z = \{ 0,1,2,3 \} $ and let $\zeta := \exp (2\pi i/8)$
be a primitive eighth root of unity in $\C$, with $i=\zeta ^2$.
Then
$${\cal C}(\rho _a) = \langle
\left(\begin{array}{cccc}
1 & 0 & 0 & 0 \\
0 & 0 & 0 & 1 \\
0 & 0 & 1 & 0 \\
0 & 1 & 0 & 0 \end{array} \right),~ \frac{1}{2}
\left(\begin{array}{rrrr}
1 & 1 & 1 & 1 \\
1 & i & -1 & -i \\
1 & -1 & 1 & -1 \\
1 & -i & -1 & i \end{array} \right), ~\diag(1,\zeta, -1,\zeta ) \rangle ~,$$
of order $|{\cal C}(\rho _a) | = 1536, $
and
$${\cal C}(\rho _b) = \langle {\cal C}(\rho _a), ~\diag (1,i,-1,-i ) \rangle ~,$$
of order $|{\cal C}(\rho _b) | = 6144 $.
Furthermore, ${\cal C}(\rho _a) \cong Z_8 \times (Z_2 \times Z_2)\,.\, SL_2(\Z/4\Z )$ and
${\cal C}(\rho _b) \cong Z_8 \times (Z_4 \times Z_4)\,.\, SL_2(\Z/4\Z )$.
The Molien series are
$$\mbox{Molien} ({\cal C}(\rho _a) ) =
\frac{(1+t^8) (1+t^{16})^2}{(1-t^8)^3(1-t^{24})} $$
and
$$\mbox{Molien} ({\cal C}(\rho _b) ) =
\frac{(1+t^{16})( 1+ t^{32})}{(1-t^8)^2(1-t^{16})(1-t^{24})} $$
(see \cite{RaSl98}).
We will interpret these rings in terms of complete weight enumerators of
self-dual isotropic codes.
If $p$ is a prime power $p\equiv \pm 1 \pmod{8} $, let
$\tilde{QR}(p) $ denote the extended quadratic residue code
of length $p+1$ over $\Z /4\Z$ (see for example \cite{BCS95}).
With the correct definition of extension,
$\tilde{QR}(p) $ contains the all-ones vector.
Let
$d_8$, $c_{16}$, $d_{16}$ be the codes of Type $\rho _b$
(see \cite{Fiel02}) with generator matrices
$$ d_8 =\begin{bmatrix}
 1 3 1 0 0 1 0 2 \\
 1 3 0 1 0 2 1 0 \\
 1 3 0 0 1 0 2 1 \\
 2 2 0 0 0 0 0 0 \\
 2 0 2 2 2 0 0 0 \end{bmatrix} \,, \ \
c_{16}= \begin{bmatrix}
1111111111111111 \\
1011111100001000 \\
1101001111000100 \\
1110101010100010 \\
0000111111100001 \\
0000020000022002 \\
0000002000022222 \\
0000000200002202 \\
0000000020000222 \\
0000000002020202 \\
0000000000220022 \end{bmatrix} \,, \ \ \
d_{16} = \begin{bmatrix}
 1 1 1 1 1 1 1 1 1 1 1 1 1 1 1 1 \\
 1 1 1 0 0 0 0 0 2 3 0 0 0 0 0 0 \\
 1 1 0 1 0 0 0 0 0 2 3 0 0 0 0 0 \\
 1 1 0 0 1 0 0 0 2 2 2 3 0 0 0 0 \\
 1 1 0 0 0 1 0 0 0 2 2 2 3 0 0 0 \\
 1 1 0 0 0 0 1 0 2 2 2 2 2 3 0 0 \\
 1 1 0 0 0 0 0 1 0 2 2 2 2 2 3 0 \\
 1 0 1 1 1 1 1 1 0 2 2 2 2 2 2 1 \\
\end{bmatrix} \, .
$$
Let $p_1:=\cwe (\tilde{QR}(7) ) $,
$p_2:=\cwe (d_8 ) $,
$p_3:=\cwe (c_{16} ) $,
$p_4:=\cwe (d_{16} ) $,
$p_5:=\cwe (\tilde{QR}(23) ) $.
Then
$$\Inv ({\cal C}(\rho _b)) = \C[p_1,p_2,p_3,p_5] (1+p_4)(1+p_6) ~,$$
where $p_6$ is the weight enumerator of a certain
code $e_{32}$ of Type $\rho _b$ and length 32.
We have an explicit description of a code $e_{16}$
that works, but it has no structure and we do not give it here.  It would
be nice to have a better example.

To find additional generators for the invariant ring of ${\cal C}(\rho _a)$,
let $e_8'$ and $d_8'$ be the codes of Type $\rho _a$ obtained from
$\tilde{QR}(7)$ resp. $d_8$ by multiplying one column by $3 \in \Z/4\Z $,
and let $p_{1a} := \cwe (e_8')$ and $p_{2a} := \cwe (d_8') $ be their complete
weight enumerators.
We also need three further codes of Type $\rho_a$:
two of length $24$ and one of length $32$, with
complete weight enumerators $f_1, f_2, f_3$ respectively.
We have examples of such codes, but again they have no structure
and we do not give them here.
With this notation we find that
 $$\Inv ({\cal C}(\rho _b)) = \C[p_1,p_2,p_{1a},p_5]
(1+p_{2a}+p_3+p_4+f_1 + f_2 + f_3+p_{2a}f_3) ~.$$

\section{Self-dual codes over $\F_{q^2} + \F_{q^2} \, u $}\label{Sec7}

In this section we study self-dual codes over the ring
$R=\F_{q^2} + \F_{q^2} \, u$,
for $q=p^f$,
where $u^2 = 0$ and $ua = a^q u$ for all $a \in \F _{q^2}$.
These are certainly ``non-classical'' codes.

We have $R\cong {\cal M}/p{\cal M}$,
where ${\cal M}$ is the maximal order in the
quaternion division algebra over the unramified extension of degree $f$
of the $p$-adic numbers.
The most important case is $q = 2$.
In this special case,
self-dual codes have been studied by Bachoc
\cite{Bach97} in connection with the construction of
interesting modular lattices, and Gaborit \cite{Gab96}
has found a mass formula.

To construct self-dual codes, we define an
$R$-valued Hermitian form $R^N \times R^N \rightarrow R$ by
$(x,y) := \sum _{i=1}^N x_i \overline{y_i} $, where
$\overline{\phantom{a}} : R \rightarrow R $ is the involution defined by
$\overline{a+bu} := a^q - bu $.
Then $$(a'+b'u)\overline{(a+bu)} =a'a^{q} + (ab'-ba') u ~,$$
for all $a,b,a',b' \in \F_{q^2}$.
A code $C\leq R^N$ is self-dual if
$C = C^{\perp } := \{ v \in R^N \mid (v,c) = 0 \mbox{~for~all~} c \in C \} $.

To express this self-duality in our language of Types, we need
a form ring.
Let $\beta :R\times R \rightarrow \frac{1}{p} \Z /\Z $ be the
bilinear form defined
by $$\beta (a'+b'u,a+bu):= \frac{1}{p} \Tr ( ab'-a'b) ~,$$
where $\Tr $ denotes the trace from $\F_{q^2}$ to $\F _p \cong \Z/p\Z$.
Let $M := \beta (1\otimes R)$, where the right action of
$R\otimes R$ on $M$ is left multiplication on the arguments:
$$ m ((r+su)\otimes (r' +s'u )) (a+bu,a'+b'u) = m((r+su)(a+bu),(r'+s'u)(a'+b'u)) ~,$$
for $m \in M$ and where $r+su, r'+s'u, a+bu, a'+b'u$ are elements of $R$.
Let $\psi : R_R \rightarrow M_{1\otimes R}$ be the $R$-module isomorphism
defined by $\psi(1) := \beta $.
The involution $J$ induced by $\beta $
is given by $(r+su)^J = r-s^qu $, and
$\epsilon = -1$ (since $\beta $ is skew-symmetric).

Define the $\frac{1}{p}\Z/\Z$-valued quadratic form
$\phi _0:R \rightarrow \frac{1}{p}\Z/\Z$ by $$ \phi _0(a+bu) := \frac{1}{p} \Tr _{\F_q/\F_p} (a a^{q}) ~,$$
and let $\Phi : = \{ \phi _0 [r] \mid r\in R \}$.
The mapping $\leftBra \rightBra : M \rightarrow \Phi $ is the obvious diagonal evaluation
$$\leftBra \beta (1\otimes (r+su)) \rightBra (a+bu) :=
\beta (1\otimes (r+su) ) (a+bu,a+bu) = -\Tr (s a a^{q}) ~.$$
Since $\leftBra \rightBra $ is surjective, this defines a unique mapping
$\lambda : \Phi \rightarrow M$, satisfying
$$\lambda \leftBra m \rightBra = m + \tau (m) ~,$$
for all $m \in M$.
To find $\lambda (\phi _0)$,
we choose $\omega \in \F_{q^2} \setminus
\F _q$ with $\omega +\omega ^q = -1$. Then
$\phi _0 = \leftBra \beta (1\otimes \omega u) \rightBra $,
$$\lambda ( \phi _0 ) = m_{0} (1\otimes \omega u ) + \tau (\beta (1\otimes \omega u)) = \beta (1\otimes u) ~,$$
and therefore $\psi ^{-1} (\lambda ( \phi _0)) = u $.
Identifying $\F_p$ with $\frac{1}{p} \Z / \Z $, we obtain
$\Q /\Z $-valued quadratic and bilinear forms.
This defines a form ring $\rho = (R,V,\beta, \Phi )$, and
the self-dual codes $C \leq R^N$ defined above are precisely
the codes of Type $\rho $.

The hyperbolic co-unitary group ${\cal U}(R,\Phi )$ contains a
normal subgroup $N$
%$$
%N:={\cal U}((u),\Phi ) := \{ (A,\phi ) \in {\cal U}(R,\Phi ) \mid
%A\equiv I_2 \pmod{u} \}
%$$
for which the quotient is a
subgroup of ${\cal U}(\F_{q^2}, \{ 0\} )$.
In fact, since $R = \F_{q^2} + \F_{q^2}u$, 
${\cal U}(R,\Phi )$ has a subgroup $H \cong O_2^+(\F_{q^2} )$
consisting of the elements
$$\{ (\begin{pmatrix} a & b \\ c & d \end{pmatrix}, \begin{pmatrix}
0 & \psi(cb) \\ \psi (da -1) & 0 \end{pmatrix} )
\mid a,b,c,d\in \F_{q^2}, ca=db=0, cb+da=1 \} ~.$$
$H$ is isomorphic to ${\cal U}(\F _{q^2},\{ 0 \})$
and is a complement to the normal
subgroup ${\cal U}((u), \Phi )$ given by
$$\{ (\begin{pmatrix} 1+au & bu \\ cu & 1+du \end{pmatrix}, \diag (\phi _0 [c'], \phi _0 [b'] )) \mid
c,b\in \F_q, a,d \in \F_{q^2}, a = d^q \} ~,$$
which is isomorphic to $\F_{q^2} \oplus \F_q \oplus \F_q $.
Therefore
$${\cal U}(R,\Phi ) \cong (\F _{q^2} \oplus \F_q \oplus \F_q ) \rtimes
O_2^+(\F_{q^2}) ~.$$

\begin{eg}
Let $q:=2$, $R= \F_4 + \F_4 u$.
Then
the hyperbolic co-unitary group ${\cal U} (R,\Phi ) $ is
generated by
$$g_1:= ( \begin{pmatrix} \omega & 0 \\ 0 & \omega ^2 \end{pmatrix}, 0 ), \
g_2:= ( \begin{pmatrix} 1+u & 0 \\ 0 & 1+u \end{pmatrix}, 0 ), $$ $$
g_3:= ( \begin{pmatrix} 1 & u \\ 0 & 1 \end{pmatrix},
 \begin{pmatrix} 0 & 0 \\ & \phi _0 \end{pmatrix}), \
h:= ( \begin{pmatrix} 0 & 1 \\ 1 & 0 \end{pmatrix},
 \begin{pmatrix} 0 & 1 \\ & 0 \end{pmatrix}) .$$
Since $\lambda $ is injective, ${\cal U}(R,\Phi )$ is isomorphic
to its image under the projection of the first component.
This image contains a normal subgroup $N \cong \F_4 + \F_2 + \F_2 $
generated by
$$\begin{pmatrix} 1+\omega u & 0 \\ 0 & 1+\omega ^2 u \end{pmatrix}, \
\begin{pmatrix} 1+ u & 0 \\ 0 & 1+ u \end{pmatrix}, \
\begin{pmatrix} 1 & u \\ 0 & 1 \end{pmatrix}, \
\begin{pmatrix} 1 & 0 \\ u & 1 \end{pmatrix} ~.$$
The quotient group is isomorphic to $S_3$, generated by the matrices
$$\begin{pmatrix} \omega & 0 \\ 0 & \omega ^2 \end{pmatrix}, \ \
\begin{pmatrix} 0 & 1 \\ 1 & 0 \end{pmatrix} .$$
Hence $${\cal U}(R,\Phi) \cong ((C_2)^2 \times (C_2)^2) \rtimes S_3 $$
where $S_3$ acts faithfully on one copy $(C_2)^2$ and
with kernel $C_3$ on the other copy.

The Molien series of ${\cal C}(\rho )$ is
$$\frac{f(t)}{(1-t^2)^5 (1-t^3) (1-t^4)^6 (1-t^6)^4} ~,$$
where
\begin{eqnarray*}
& f(t) := & 1 + t + 4t^2 + 3t^3 + 53t^4 + 104t^5 + 458t^6 + 858t^7 + 2474t^8 + 4839t^9 \\
&&{} ~+ 10667t^{10 } + 19018t^{11 } + 34193t^{12 } + 55481t^{13 } + 86078t^{14 } \\
&&{} ~+ 125990t^{15 } + 173466t^{16 } + 230402t^{17 } + 287430t^{18 } + 346462t^{19 } \\
&&{} ~+ 393648t^{20 } + 431930t^{21 } + 450648t^{22 } + 450648t^{23 } + 431930t^{24 } \\
&&{} ~+ 393648t^{25 } + 346462t^{26 } + 287430t^{27 } + 230402t^{28 } + 173466t^{29 } \\
&&{} ~+ 125990t^{30 } + 86078t^{31 } + 55481t^{32 } + 34193t^{33 } + 19018t^{34 } \\
&&{} ~+ 10667t^{35 } + 4839t^{36} + 2474t^{37 } + 858t^{38 } + 458t^{39 } + 104t^{40 } \\
&&{} ~+ 53t^{41 } + 3t^{42 } + 4t^{43 } + t^{44 } + t^{45 } ~.
\end{eqnarray*}

Various interesting symmetrizations are possible:
\begin{itemize}
\item[a)] Symmetrizing by the action of $\F _4^* = \langle \omega \rangle$
yields a matrix group ${\cal C}_{(7)}(\rho ) \cong (C_2)^2 \rtimes S_3$
(a non-faithful $S_3$-action) of degree 7,
with Molien series
$$\frac{1+t^2+9t^4+21t^6+41t^8+43t^{10}+43t^{12}+23t^{14}+10t^{16}}{(1-t)(1-t^2)(1-t^4)^3(1-t^6)^2} ~.$$
\item[b)] The unit group $R^*$ has three orbits on $R$
namely $\{0\}, R^*, uR^* $.
Symmetrizing by $R^*$ gives a matrix group ${\cal C}_{(3)}(\rho ) \cong D_8$,
for which the Molien series and invariant ring were described by
Bachoc \cite[Theorem 4.4]{Bach97}.
\end{itemize}
\end{eg}

\end{document}